\date{today}

\documentclass[10pt, letterpaper]{amsart}

\usepackage{bm}
\usepackage{latexsym, amsmath, amstext, amssymb, amsfonts, amscd, bm, array, multirow, amsbsy, mathrsfs}
\usepackage{amsthm}
\usepackage{t1enc}
\usepackage[mathscr]{eucal}
\usepackage{indentfirst}
\usepackage{pb-diagram}
\usepackage{graphicx}
\usepackage{fancyhdr}
\usepackage{fancybox}
\usepackage{enumerate}
\usepackage{color}
\usepackage[all]{xy}
\usepackage{hyperref}
\usepackage{tikz}
\usepackage{xparse}
\hypersetup{colorlinks=false,pdfborderstyle={/S/U/W 0}}
\usetikzlibrary{matrix}

\newcommand{\arxiv}[1]{{\tt
		\href{http://www.arXiv.org/abs/#1}{arXiv:#1}}}

\usepackage{url}
\usepackage[sort&compress,comma]{natbib}
\bibpunct{[}{]}{,}{n}{}{,}
\theoremstyle{plain}
\newtheorem{thm}{Theorem}[section]

\newtheorem{cor}[thm]{Corollary}

\theoremstyle{definition}
\newtheorem{definition}[thm]{Definition}

\theoremstyle{remark}
\newtheorem{remark}[thm]{Remark}

\newcommand{\Hom}{{\rm Hom}}

\def\ad{\mathrm{ad}}

\newcommand{\bDelta}{\boldsymbol{\Delta}}

\newcommand{\bl}{\begin{Lemma}}
	\newcommand{\el}{\end{Lemma}}
\newcommand{\bt}{\begin{Theorem}}
	\newcommand{\et}{\end{Theorem}}
\newcommand{\bd}{\begin{Definition}}
	\newcommand{\ed}{\end{Definition}}
\newcommand{\End}{\mathrm{End}}
\newcommand{\Aut}{\mathrm{Aut}}

\newcommand{\ES}{\mathrm{ES}}
\newcommand{\DS}{\mathrm{DS}}

\newcommand{\Symp}{\mathrm{Symp}}

\newcommand{\eqdef}{\stackrel{{\rm def.}}{=}}

\DeclareFontFamily{U}{rsf}{}
\DeclareFontShape{U}{rsf}{m}{n}{<5> <6> rsfs5 <7> <8> <9> rsfs7 <10-> rsfs10}{}
\DeclareMathAlphabet\Scr{U}{rsf}{m}{n}

\def\cU{\mathcal{U}}

\def\Z{\mathbb{Z}}
\def\C{\mathbb{C}}
\def\R{\mathbb{R}}

\def\rk{{\rm rk}}

\def\dd{\mathrm{d}}

\def\ad{\mathrm{ad}}

\def\bcD{\boldsymbol{\mathcal{D}}}

\def\bt{\mathbf{t}}

\def\U{\mathrm{U}}

\newcommand{\be}{\begin{equation*}}
\newcommand{\ee}{\end{equation*}}
\newcommand{\ben}{\begin{equation}}
\newcommand{\een}{\end{equation}}
\newcommand{\beqa}{\begin{eqnarray*}}
	\newcommand{\eeqa}{\end{eqnarray*}}
\newcommand{\beqan}{\begin{eqnarray}}
\newcommand{\eeqan}{\end{eqnarray}}

\newcommand{\id}{\mathrm{id}}
\newcommand{\Tr}{\mathrm{Tr}}

\def\Diff{\mathrm{Diff}}

\def\cC{{\mathcal C}}

\def\loslash{\bm{\oslash}}

\newcommand{\cM}{\mathcal{M}}

\def\Hol{\mathrm{Hol}}

\def\cD{\mathcal{D}}

\def\cE{\check{E}}
\def\cX{\mathcal{X}}

\def\cP{\mathcal{P}}

\def\cG{\mathcal{G}}
\def\cT{\mathcal{T}}

\def\cC{\mathcal{C}}
\def\cH{\mathcal{H}}
\def\bXi{{\boldsymbol{\Xi}}}

\def\Sp{\mathrm{Sp}}
\def\G_2{\mathrm{G_2}}

\def\cS{\mathcal{S}}

\def\cV{\mathcal{V}}

\def\G{\mathrm{G}}

\def\cE{\mathcal{E}}

\def\iOmega{\mathit{\Omega}}

\def\cX{\mathcal{X}}

\def\grad{\mathrm{grad}}

\def\bd{\boldsymbol{\dd}}
\def\bast{\boldsymbol{\ast}}

\def\Met{\mathrm{Met}}
\def\Iso{\mathrm{Iso}}

\def\fJ{\mathfrak{J}}

\def\i{\mathbf{i}}
\def\sc{\mathrm{sc}}

\def\Sol{\mathrm{Sol}}
\def\Conf{\mathrm{Conf}}
\def\ub{\mathrm{ub}}

\def\bwedge{\boldsymbol{\wedge}}

\def\cl{\mathrm{cl}}

\def\bAd{\mathbf{Ad}}

\def\bcX{\boldsymbol{\cX}}

\begin{document}

\title{The global formulation of generalized Einstein-Scalar-Maxwell theories}

\author[C. Lazaroiu]{C. Lazaroiu}
\address{Center for Geometry and Physics, Institute for Basic
	Science, Pohang, Republic of Korea}
\email{calin@ibs.re.kr}

\author[C. S. Shahbazi]{C. S. Shahbazi}
\address{Institut de Physique Th\'eorique, CEA-Saclay, France}
\email{carlos.shabazi-alonso@cea.fr}


\begin{abstract}
We summarize the global geometric formulation of Einstein-Scalar-Maxwell theories twisted by flat symplectic vector bundle which encodes the duality structure of the theory. We describe the scalar-electromagnetic symmetry group of such models, which consists of flat unbased symplectic automorphisms of the flat symplectic vector bundle lifting those isometries of the scalar manifold which preserve the scalar potential. The Dirac quantization condition for such models involves a local system of integral symplectic spaces, giving rise to
a bundle of polarized Abelian varieties equipped with a symplectic
flat connection, which is defined over the scalar manifold of the
theory. Generalized Einstein-Scalar-Maxwell models arise as the bosonic sector of the effective theory of string/M-theory compactifications to four-dimensions, and they are characterized by having non-trivial solutions of ``U-fold'' type.
\end{abstract}

\maketitle

\setcounter{tocdepth}{1} 
\tableofcontents

\section{Introduction}

Supergravity theories \cite{Ortin,FreedmanProeyen} are supersymmetric
theories of gravity which extend general relativity and gauge theory
and arise in the low energy limit of string/M-theory and of their
compactifications. It is known that the construction of such theories
involves interesting structures, such as K\"ahler-Hodge and special
K\"ahler manifolds, symmetric spaces, exceptional Lie groups,
generalized complex structures, differential cohomology and
differential K-theory etc. However, the global formulation of
these theories is not yet fully understood. This note is part of a
larger project (see also
\cite{GeometricUfolds,Lipschitz,lip}) aimed at obtaining the
full mathematical formulation of supergravity theories (in the
generality required by their relation to string theory) and at studying
the global geometry of their solutions.

Supergravity theories are classical theories of gravity coupled to
matter, formulated using systems of ``fields'' defined on a manifold
$M$ of the appropriate dimension and subject to certain partial
differential equations, known as the ``equations of motion''. An
unambiguous formulation of such theories requires that one specifies
the global nature of the fields and of the partial differential
operators arising in the equations of motion. Currently, however, the
supergravity literature gives only {\em local}
descriptions\footnote{Descriptions which are valid only if one
restricts all fields to sufficiently small open subsets of $M$.} of
the fields and of these differential operators.  The {\em
globalization problem} is the problem of giving globally-unambiguous
mathematical definitions of such theories which reduce locally to the
local description found in the supergravity literature. The solution
of this problem is non-unique since there can be many global
definitions of ``fields'' subject to globally-defined partial
differential equations which reduce to a given local description.

Since supergravity theories are supersymmetric, they require spinors
for their formulation. In this note, we simplify the globalization problem by
ignoring the spinor field content and the supersymmetry conditions,
thus considering only the so-called
\emph{universal bosonic sector}. This sector arises in any
supergravity theory, though it is subject to increasingly stringent
supplementary constraints (not discussed in this paper) as the number
of supersymmetries present in the theory increases. In addition, we
focus exclusively on the case when $M$ is a four-manifold.  

In four dimensions, the universal bosonic sector is the so-called
Einstein-Scalar-Maxwell (ESM) model defined on a four-manifold $M$,
which involves gravity (modeled globally by a Lorentzian metric on
$M$), a finite number of real scalar fields (modeled globally by a
smooth map from $M$ to a manifold $\cM$ of arbitrary dimension) and a
finite number of Abelian gauge fields, whose field strengths can be
modeled {\em locally} as 2-forms defined on $M$. While the local form
of ESM theories is well-known, their precise global formulation was
systematically addressed only recently \cite{sigma}. It
turns out that the naive globalization of the local formulation
fails to capture the classical limit of certain string theory
backgrounds known as ``U-folds'' and hence is insufficient for the
application of such models to string theory. The geometric description
of the classical limit of U-fold backgrounds \cite{GeometricUfolds}
requires that one globalizes ESM models by including a ``twist'' of
the Abelian gauge field sector through the (pull-back of) a flat
symplectic vector bundle defined on $\cM$. This
produces so-called {\em generalized ESM models}, which are {\em locally}
indistinguishable from the naive globalization but have a rather
different global behavior. The naive globalization corresponds to
using a {\em trivial} flat symplectic vector bundle on $\cM$.

The global mathematical formulation of generalized ESM models given in 
\cite{sigma} is summarized in this note. We follow the
notations and conventions of loc. cit.; in particular, all manifolds
considered are smooth and connected and all bundles and maps
considered are smooth. In this note, a Lorentzian metric is a smooth
metric of signature $(3,1)$ defined on a four-manifold.

\section{Generalized Einstein-Scalar-Maxwell theories} 

\subsection{Scalar structures and related notions}

\begin{definition} 
A {\em scalar structure} is a triplet $\Sigma=(\cM,\cG,\Phi)$, where
$(\cM,\cG)$ is a Riemannian manifold (called the {\em scalar
manifold}) and $\Phi\in \cC^\infty(\cM,\R)$ is a smooth real-valued
function defined on $\cM$ (called the {\em scalar potential}).
\end{definition}

\noindent Let $\Sigma=(\cM,\cG,\Phi)$ be a scalar structure. 
Let $M$ be an oriented four-dimensional smooth manifold 
(which need not be compact).

\begin{definition}
The {\em modified density} of a smooth map
$\varphi\in \cC^\infty(M,\cM)$ relative to a Lorentzian metric
$g\in \Met_{3,1}(M)$ and to the scalar structure $\Sigma$ is the
following smooth real-valued map defined on $M$:
\ben 
e_\Sigma(g,\varphi)\eqdef \frac{1}{2}\Tr_g \varphi^\ast(\cG)+\Phi^\varphi\in \cC^\infty(M,\R)~~,
\een
where $\Phi^{\varphi} \eqdef \Phi\circ\varphi$ and $\Tr_g$ denotes
trace taken with respect to $g$.
\end{definition}

\begin{definition}
The {\em modified tension field} of a smooth map $\varphi\in
\cC^\infty(M,\cM)$ relative to the Lorentzian metric $g\in \Met_{3,1}(M)$ and
to the scalar structure $\Sigma$ is the section of the pulled-back
bundle $(T\cM)^\varphi$ defined through:
\begin{equation}
\theta_\Sigma(g,\varphi)\eqdef \theta_\cG(g,\varphi)-(\grad_\cG \Phi)^\varphi\in \Gamma(M, (T\cM)^\varphi)\, .
\end{equation}
Here $\grad_\cG \Phi\in \cX(\cM)$ is the gradient vector field of
$\Phi$ with respect to $\cG$ and $\theta_\cG(g,\varphi)$ is the
tension field of $\varphi$ relative to $g$ and $\cG$ \cite{BairdWood}:
\begin{equation}
\theta_{\cG}(g,\varphi) \eqdef \Tr_{g} \nabla \widetilde{\dd\varphi} \in \Omega^{0}(M, (T\cM)^\varphi)\, ,
\end{equation}
where $\widetilde{\dd\varphi}\in \Omega^{1}(M,(T\cM)^{\varphi})$ denotes
the $(T\cM)^{\varphi}$-valued one-form associated to the differential
$\dd\varphi\colon TM\to T\cM$ and $\nabla$ is the connection induced on
$(T\cM)^\varphi$ by the Levi-Civita connections of $g$ and $\cG$.
\end{definition}

\subsection{Duality structures}

\noindent Let $N$ be a manifold.

\begin{definition}
A {\em duality structure} on $N$ is a flat symplectic vector bundle
$\Delta=(\cS,D,\omega)$ defined over $N$, where $\omega$ denotes the
symplectic pairing on the vector bundle $\cS$ and $D$ denotes the
$\omega$-compatible flat connection on $\cS$.
\end{definition}

\begin{definition}
Let $\Delta_{i} = (\cS_{i},D_{i},\omega_{i})$ with $i=1,2$ be two
duality structures defined on $N$. A {\em morphism of duality
structures} from $\Delta_{1}$ to $\Delta_{2}$ is a based morphism of
vector bundles $f\in \Hom(\cS_1,\cS_2)$ such that $\omega_2(f\otimes
f)=\omega_1$ and such that $D_2\circ f=(\id_{\Omega^1(N)}\otimes
f)\circ D_1$.
\end{definition}

\noindent With this notion of morphism, duality structures on $N$ form a
category which we denote by $\mathrm{DS}(N)$. Let
$\Delta=(\cS,D,\omega)$ be a duality structure defined on $N$ such
that $\rk\,\cS=2n$.  Let $\Symp$ denote the category of
finite-dimensional symplectic vector spaces over $\R$ and linear
symplectic morphisms. Let $\Symp^{\times}$ denote the unit groupoid of
this category and $\Pi_1(N)$ denote the fundamental groupoid of
$N$. Let $T_\gamma^\Delta$ denote the parallel transport of $D$ along
a path $\gamma:[0,1]\rightarrow N$.

\begin{definition}
The {\em parallel transport functor} of $\Delta$ is the functor
$T_{\Delta}:\Pi_1(N)\rightarrow \Symp^\times$ which associates to any
point $x\in N$ the symplectic vector space
$T_\Delta(x)=(\cS_x,\omega_x)$ and to any homotopy class $c\in
\Pi_1(N)$ with fixed initial point $x$ and fixed final point $y$ the
invertible symplectic morphism
$T_\Delta(c)=T_\gamma^\Delta:(\cS_x,\omega_x)\stackrel{\sim}{\rightarrow}
(\cS_y,\omega_y)$, where $\gamma\in \cP(N)$ is any path which
represents the class $c$.
\end{definition}

\noindent Notice that $T_\Delta$ can be viewed as a $\Symp^\times$-valued local system
defined on $N$. The map which takes $\Delta$ into $T_\Delta$ is an
equivalence between the category $\mathrm{DS}(N)$ and the functor
category $[\Pi_1(N),\Symp^\times]$. This implies that duality
structures on $N$ are classified up to isomorphism by the symplectic
character variety:
\begin{equation}
C_{\pi_1(N)}(\Sp(2n,\R))\eqdef \Hom(\pi_1(N),\Sp(2n,\R))/\Sp(2n,\R)~~.
\end{equation}

\begin{definition}
A {\em duality frame} of $\Delta$ is a $D$-flat symplectic frame
$\cE\eqdef (e_1,\ldots e_n,f_1,\ldots , f_n)$ of $(\cS,\omega)$ defined on an
open subset $\cU\subset N$.
\end{definition}

\begin{definition}
The duality structure $\Delta$ is called {\em trivial} if it is
trivial as a flat symplectic vector bundle.
\end{definition}

\begin{remark}
A duality structure is trivial iff it admits a globally-defined duality frame. 
If $N$ is simply connected, then any duality structure on $N$ is trivial. 
\end{remark}

\subsection{Electromagnetic structures}

\noindent Let $N$ be a manifold.

\begin{definition}
An {\em electromagnetic structure} defined on $N$ is a quadruplet
$\Xi\eqdef (\cS,D,J,\omega)$, where $(\cS,D,\omega)$ is a duality
structure defined on $N$ and $J$ is a taming of the symplectic vector
bundle $(\cS,\omega)$.
\end{definition}

\begin{remark}
Notice that we do not require $J$ to be compatible with $D$. Together
with $\omega$, $J$ defines an Euclidean scalar product $Q$ on $\cS$ given by 
$Q(\cdot,\cdot) \eqdef \omega(J\cdot,\cdot)$.
\end{remark}

\begin{definition}
Let $\Xi_1=(\cS_1,D_1,J_1,\omega_1)$ and
$\Xi_2=(\cS_1,D_1,J_1,\omega_1)$ be two electromagnetic structures
defined on $N$. A {\em morphism of electromagnetic structures} from
$\Xi_1$ to $\Xi_2$ is a morphism of duality structures
$f:(\cS_1,D_1,\omega_1)\rightarrow (\cS_2,D_2,\omega_2)$ such that
$J_2\circ f =f\circ J_1$.
\end{definition}

\noindent With this definition of morphism, electromagnetic structures defined 
on $N$ form a category which we denote by $\ES(N)$. This fibers over
the category of duality structures $\DS(N)$; the fiber at a duality
structure $\Delta=(\cS,D,\omega)$ can be identified with the set
$\fJ_+(\cS,\omega)$ of tamings of $(\cS,\omega)$, which is a
contractible topological space. Accordingly, the set of isomorphism
classes of $\ES(N)$ fibers over the disjoint union of character
varieties $\sqcup_{n\geq 0}C_{\pi_1(N)}(\Sp(2n,\R))$. Let
$\Xi=(\cS,D,J,\omega)$ be an electromagnetic structure defined on $N$
and $h=Q+\i\omega$ be the Hermitian scalar product defined by $\omega$
and $J$ on $\cS$.

\begin{definition}
The {\em fundamental form} of $\Xi$ is the $End(\cS)$-valued one-form
on $N$ defined through:
\be
\Theta_\Xi \eqdef D^\ad(J)\eqdef D\circ J-J\circ D\in \Omega^1(N,End(\cS))~~.
\ee
The electromagnetic structure $\Xi$ is called \emph{unitary} if
$\Theta_\Xi=0$, i.e. if $J$ is parallel with respect to $D$.
\end{definition}

\noindent If $\Xi$ is unitary, then $D$ is a unitary connection 
on the Hermitian vector bundle $(\cS,J,h)$. In this case, we have
$\Hol_D^x\subset \U(\cS_x,J_x,h_x)$ for all $x\in N$. The {\em
category of unitary electromagnetic structures} defined on $N$ is the
full sub-category of $\ES(N)$ whose objects are the unitary
electromagnetic structures.  This is equivalent with the category of
Hermitian vector bundles defined on $N$ and endowed with a flat
$\C$-linear Hermitian connection. In particular, isomorphism classes
of {\em unitary} electromagnetic structures are in bijection with the
points of the character variety:
\be
C_{\pi_1(N)}(\U(n))\eqdef \Hom(\pi_1(N),\U(n))/\U(n)~~,
\ee
where $\U(n)$ acts by conjugation. 

\subsection{Scalar-duality and scalar-electromagnetic structures}

\begin{definition} 
A {\em scalar-duality structure} is an ordered system $(\Sigma,\Xi)$, where:
\begin{enumerate}[1.]
\itemsep 0.0em
\item $\Sigma=(\cM,\cG,\Phi)$ is a scalar structure
\item $\Xi=(\cS,D,\omega)$ is a duality structure defined on
  $\cM$.
\end{enumerate}
A {\em scalar-electromagnetic structure} is an ordered system
$\cD=(\Sigma,\Xi)$, where:
\begin{enumerate}[1.]
\itemsep 0.0em
\item $\Sigma=(\cM,\cG,\Phi)$ is a scalar structure
\item $\Xi=(\cS,D,J,\omega)$ is an electromagnetic structure defined on
  $\cM$.
\end{enumerate}
In this case, the system $\cD_0\eqdef (\Sigma,\Xi_0)$ is called the
  underlying scalar-duality structure, where $\Xi_0\eqdef (\cS,D,\omega)$ is
  the duality structure underlying $\Xi$.
\end{definition}

\noindent Let $\cD$ be a scalar-electromagnetic structure as in the definition. 

\begin{definition}
The {\em fundamental field} of the scalar-electromagnetic structure
$\cD$ is defined through:
\be
\Psi_\cD \eqdef (\sharp_\cG\otimes \id_{End(\cS)})(\Theta_\Xi) \in \Gamma(\cM, T\cM\otimes End(\cS))\, . 
\ee
\end{definition}

\subsection{Pulled-back electromagnetic structures}

\noindent Let $\cD=(\Sigma,\Xi)$ be a scalar-electromagnetic structure with
underlying scalar structure $\Sigma=(\cM,\cG,\Phi)$ and underlying
electromagnetic structure $\Xi=(\cS,D,J, \omega)$.  Let $M$ be a
four-manifold and $\varphi\in \cC^\infty(M,\cM)$ be a smooth map from
$M$ to $\cM$.

\begin{definition}
The $\varphi$-pullback of the electromagnetic structure $\Xi$ defined
on $\cM$ is the electromagnetic structure
$\Xi^\varphi\eqdef(\cS^\varphi, D^\varphi,J^\varphi,\omega^\varphi)$
defined on $M$.
\end{definition}

\noindent The Hodge operator $\ast_g:\wedge T^\ast M\rightarrow \wedge
T^\ast M$ of $(M,g)$ induces the endomorphism $\bast\eqdef
\bast_g\eqdef \ast_g\otimes \id_{\cS^{\varphi}}$ of the bundle
$\wedge_M(\cS^\varphi)\eqdef \wedge T^\ast M\otimes \cS^{\varphi}$.

\begin{definition}
The {\em twisted Hodge operator} of $\Xi^\varphi$ is the bundle
endomorphism $\star:=\star_{g,J^\varphi}\in \End(M,\wedge T^\ast
M\otimes \cS^\varphi)$ defined through:
\be
\star_{g,J^\varphi}\eqdef\ast_g\otimes J^\varphi=\bast_g\circ J^\varphi=J^\varphi\circ \bast_g~~.
\ee
\end{definition}

\noindent Let $\alpha \eqdef \oplus_{k=0}^4 (-1)^k\id_{\wedge^k T^\ast
  M}$ be the {\em main automorphism} of $\wedge T^\ast M$. We have:
\ben
\label{HodgeSquare}
\star^2=\alpha\otimes \id_{\cS^\varphi}~~.
\een
The operator $\star_{g,J^\varphi}$ preserves the sub-bundle
$\wedge^2_M(\cS^\varphi)=\wedge^2T^\ast M\otimes \cS^\varphi$, on
which it squares to plus the identity. Accordingly, we have a direct
sum decomposition:
\be
\wedge^2 T^\ast M\otimes \cS^\varphi=(\wedge^2 T^\ast M\otimes \cS^\varphi)^+\oplus (\wedge^2 T^\ast M\otimes \cS^\varphi)^-~~,
\ee 
where $(\wedge^2 T^\ast M\otimes \cS^\varphi)^\pm$ are the sub-bundles
of eigenvectors of $\star$ corresponding to the eigenvalues $\pm 1$.

\begin{definition}
An $\cS^\varphi$-valued two-form $\eta\in \Omega^2(M,\cS^\varphi)$
defined on $M$ is called {\em positively polarized} with respect to
$g$ and $J^\varphi$ if it is a section of the vector bundle $(\wedge^2
T^\ast M\otimes \cS^\varphi)^+$, which amounts to the requirement that
it satisfies the {\em positive polarization condition}:
\ben
\label{2formpol}
\star_{g,J^\varphi} \eta=\eta~~\mathrm{i.e.}~~\ast_g \eta=-J^\varphi\eta~~.
\een
\end{definition}

\noindent For any open subset $U$ of $M$, let $g_U\eqdef g|_U$,
$\varphi_U\eqdef \varphi|_U$ and let $\iOmega^{\Xi, g, \varphi}$ be the
sheaf of smooth sections of the bundle $(\wedge^2 T^\ast M\otimes
\cS^\varphi)^+$. Globally-defined and positively-polarized
$\cS^\varphi$-valued forms are the global sections of this
sheaf. Notice that $\eta\in \Omega^2(M,\cS^\varphi)$ is positively
polarized iff $\star \eta$ is.

\subsection{The mathematical formulation of generalized ESM theories} 

Let $M$ be a four-manifold and $\cD=(\Sigma,\Xi)$ be a
scalar-electromagnetic structure with underlying scalar structure
$\Sigma=(\cM,\cG,\Phi)$ and underlying electromagnetic structure
$\Xi=(\cS,D,J, \omega)$. The $\varphi$-pullback $Q^\varphi$ of the
Euclidean scalar product $Q$ induced by $\omega$ and $J$ on $\cS$ is a
Euclidean scalar product on $\cS^\varphi$. Let
$\oslash_g:\otimes^4T^\ast M\rightarrow \otimes^2 T^\ast M$ be the
bundle morphism given by $g$-contraction of the two middle indices.
This is uniquely determined by the condition:
\be
(\omega_1\otimes\omega_2)\oslash (\omega_3\otimes \omega_4)=(\omega_2,\omega_3)_g\omega_1\otimes \omega_4~~\forall \omega_1,\omega_2,\omega_3,\omega_4~~\forall \omega\in \Omega^1(\cM)~~,
\ee
where $(~,~)_g$ is the pseudo-Euclidean metric induced by $g$ on
$\wedge T^\ast M$. Viewing $\wedge^2 T^\ast M$ as the sub-bundle of
antisymmetric 2-tensors inside $\otimes^2 T^\ast M$, this restricts to
a morphism of vector bundles $\oslash_g: \wedge^2 T^\ast
M \otimes \wedge^2 T^\ast M\rightarrow \otimes^2 T^\ast M$, which we
call the {\em inner $g$-contraction of 2-forms}.

\begin{definition}
\label{def:twistedinner}
The {\em twisted inner contraction} of $\cS^\varphi$-valued 2-forms is
the unique morphism of vector bundles
$\loslash:=\loslash_{g, J,\omega,\varphi}:\wedge_M^2(\cS^\varphi)\times_M\wedge_M^2(\cS^\varphi)\rightarrow
\otimes^2(T^\ast M)$ which satisfies:
\be
(\rho_1\otimes \xi_1)\loslash (\rho_2\otimes \xi_2)= Q^\varphi(\xi_1,\xi_2)\rho_1\oslash_g\rho_2~~
\ee
for all $\rho_1,\rho_2\in \Omega^2(M)$ and all $\xi_2,\xi_2\in\Gamma(M,
\cS^\varphi)$. 
\end{definition}

\noindent Let $\Psi\eqdef \Psi_\cD\in \Gamma(\cM,T\cM\otimes End(\cS))$ be the
fundamental field of $\cD$ and let $\Psi^\varphi\in
\Gamma(M,(T\cM)^\varphi\otimes End(\cS^\varphi))$ be its pullback
through $\varphi$. Let $(~,~)$ be the pseudo-Euclidean scalar product
induced by $g$ and $Q^\varphi$ on the vector bundle
$\wedge_M(\cS^\varphi)\eqdef \wedge T^\ast M\otimes \cS^\varphi$. For
any vector bundle $\cT$ defined on $M$, we extend this trivially to a
$\cT$-valued pairing (denoted by the same symbol) between the bundles
$\cT\otimes \wedge_M(\cS^\varphi)$ and
$\wedge_M(\cS^\varphi)$. Similarly, we trivially extend the twisted
wedge product $\bwedge_\omega$ defined in Appendix C of
reference \cite{sigma} to a $\cT\otimes \wedge T^\ast M$-valued
pairing (denoted by the same symbol) between the bundles
$\cT \otimes \wedge_M(\cS^\varphi)$ and $\wedge_M(\cS^\varphi)$.

\begin{definition}
The {\em sheaf of ESM configurations} $\Conf_\cD$ determined by $\cD$ is
the sheaf of sets defined on $M$ through:
\be
\Conf_\cD(U)\eqdef \{(g,\varphi,\cV)|g\in \Met_{3,1}(U),\varphi\in \cC^\infty(U,\cM),\cV\in \iOmega^{\Xi,g,\varphi}(U)\}~~
\ee
for all open subsets $U\subset M$, with the obvious restriction
maps. An element $(g,\varphi,\cD)\in \Conf_\cD(U)$ is called a {\em
local ESM configuration} of type $\cD$ defined on $U$. The {\em set of
global configurations} of type $\cD$ is the set:
\be
\Conf_\cD(M)\eqdef \{(g,\varphi,\cV)|g\in \Met_{3,1}(M),\varphi\in \cC^\infty(M,\cM),\cV\in \iOmega^{\Xi,g,\varphi}(M)\}~~.
\ee
of global sections of this sheaf. An element $(g,\varphi,\cD)\in
\Conf_\cD(M)$ is called a {\em global ESM configuration} of type
$\cD$.
\end{definition}

\begin{definition}
\label{def:generalizedESM}
The \emph{generalized ESM theory} associated to $\cD$ is defined by the
following set of partial differential equations on $M$ with unknowns
$(g,\varphi,\cV)\in \Conf_\cD(M)$:
\begin{enumerate}[1.]
\item The {\em Einstein equation}:
\ben
\label{eins}
\G(g) =\kappa\, \mathrm{T}(g,\varphi,\cV)~~,
\een
with {\em energy-momentum tensor} $\mathrm{T}_\cD$ given by:
\begin{equation}
\mathrm{T}_\cD(g,\varphi,\cV) \eqdef  g\, e_\Sigma (g,\varphi)  + 2~\cV\loslash \cV - \varphi^\ast(\cG) \, . 
\end{equation}
\item The {\em scalar equations}:
\ben
\label{sc}
\theta_\Sigma(g,\varphi) - \frac{1}{2} (\bast \cV , \Psi^\varphi\cV) = 0~~.
\een
\item The {\em twisted electromagnetic equations}:
\ben
\label{em}
\dd_{D^{\varphi}}\cV = 0~~,
\een
where
$\dd_{D^\varphi}:\Omega^k(M,\cS^\varphi)\rightarrow \Omega^{k+1}(M,\cS^\varphi)$
is the de Rham differential of $M$ twisted by the pulled-back
flat connection $D^\varphi$.
\end{enumerate}
A {\em local ESM solution} of type $\cD$ defined on $U$ is a smooth
solution $(g,\varphi,\cV)$ of these equations which is defined on
$U$. A {\em global ESM solution} of type $\cD$ is a smooth solution of
these equations which is defined on $M$. The {\em sheaf of local ESM
solutions} $\Sol_\cD$ of type $\cD$ is the sheaf of sets defined on
$M$ whose sections on an open subset $U\subset M$ is the set of all
local solutions defined on $U$.
\end{definition}

\begin{remark}
\label{remark:physical}
It is shown in \cite{sigma} that a generalized ESM model is {\em
locally} indistinguishable from an ordinary ESM model, in the sense
that the global partial differential equations \eqref{eins},
\eqref{sc} and \eqref{em} reduce locally to those used in the supergravity literature
(see for example reference \cite{Ortin}\footnote{Notice however that
we use different conventions.}) upon choosing a local flat symplectic
frame of the duality structure $\Delta=(\cS,D,\omega)$.  The
supergravity literature tacitly assumes that the local formulas
globalize trivially, which amounts to working with a {\em trivial}
duality structure; this assumption implies existence of a
globally-defined duality frame. Generalized ESM models with a
non-trivial duality structure are globally quite different from the
models used in the supergravity literature, since a non-trivial
duality structure does not admit global duality frames. Due to this
fact, global solutions of generalized ESM models afford a geometric
description of a certain type of classical U-folds, thereby realizing
the proposal of \cite{GeometricUfolds}.
\end{remark}

\subsection{Sheaves of scalar-electromagnetic configurations and
  solutions}

Let $M$ and $\cD$ be as above and fix a metric $g\in \Met_{3,1}(M)$. 

\begin{definition}
The {\em sheaf of local scalar-electromagnetic configurations}
$\Conf^g_{\cD}$ relative to $g$ is the sheaf of sets defined on $M$
whose set of sections on an open subset $U\subset M$ is defined through: 
\be
\Conf^g_{\cD}(U) \eqdef \{(\varphi,\cV)|\varphi\in \cC^\infty(U,\cM),\cV\in \iOmega^{\Xi,g,\varphi}(U)\}~~.
\ee
The {\em set of global scalar-electromagnetic configurations} relative
to $g$ is the set $\Conf^g_\cD(M)$ of global sections of this sheaf.
\end{definition}

\begin{definition}
The {\em sheaf of local scalar-electromagnetic solutions} relative to
$g$ is the sheaf of sets defined on $M$ whose set of sections
$\Sol^g_\cD(U)$ on a an open subset $U\subset M$ is defined as the set
of all solutions of the scalar and twisted electromagnetic
equations \eqref{sc} and \eqref{em} defined on $U$. The {\em set of
global scalar-electromagnetic solutions} relative to $g$ is the set
$\Sol^g_\cD(M)$ of global sections of $\Sol^g_\cD$.
\end{definition}

\noindent Since it will be of use later, we define:
\ben
\label{Confgdef}
\Conf^g_{\cD_0}(M)\eqdef \cup_{J\in \fJ_+(\cS,\omega)} \Conf^g_{(\cD_0,J)}(M)~~,
\een
where $\cD_0$ is a scalar-duality structure. 


\subsection{Electromagnetic field strengths}


\begin{definition}
An {\em electromagnetic field strength} on $M$ with respect to $\cD$
and relative to $g\in \Met_{3,1}(M)$ and to the map $\varphi\in
\cC^\infty(M,\cM)$ is an $\cS^\varphi$-valued 2-form $\cV\in
\Omega^2(M,\cS^\varphi)$ which satisfies the following two conditions:
\begin{enumerate}[1.]
\itemsep 0.0em
\item $\cV$ is positively polarized with respect to $J^\varphi$,
  i.e. we have $\star_{g,J^\varphi}\cV=\cV$.
\item $\cV$ is $\dd_{D^\varphi}$-closed, i.e.:
\ben
\label{GaugeEOM0}
\dd_{D^\varphi} \cV=0~~.
\een
\end{enumerate}
The second condition is called {\em the electromagnetic equation}. 
\end{definition}

\

\noindent For any open subset $U$ of $M$, let:
\ben
\iOmega^{\Xi,g,\varphi}_\cl(U)\eqdef \{\cV\in \iOmega^{\Xi,g, \varphi}(U) | \dd_{D^\varphi} \cV=0\}
\een
denote the set of electromagnetic field strengths defined on $U$,
which is an (infinite-dimensional) subspace of the $\R$-vector space
$\iOmega^{\Xi,g,\varphi}(U)$. This defines a {\em sheaf of
  electromagnetic field strengths} $\iOmega^{\Xi,g,\varphi}_\cl$
relative to $\varphi$ and $g$, which is a locally-constant sheaf of
$\R$-vector spaces defined on $M$.

\section{Scalar-electromagnetic dualities and symmetries} 
\label{sec:DualityTransformations}

Let $\Delta=(\cS,D,\omega)$ be a duality structure on $\cM$ and $J$ be
a taming of $(\cS,\omega)$. Let $\Xi=(\cS,D,J,\omega)$ be the
corresponding electromagnetic structure with underlying duality
structure $\Delta=(\cS,D,\omega)$.

\begin{definition}
An {\em unbased} automorphism $f\in \Aut^\ub(\cS)$ is called:
\begin{enumerate}[1.]
\itemsep 0.0em
\item A {\em symmetry of the duality structure} $\Delta$, if $f$ is symplectic with 
respect to $\omega$ and covariantly constant with respect to $D$.
\item A {\em symmetry of the electromagnetic structure} $\Xi$, if $f$ is 
complex with respect to $J$ and is a symmetry of the duality structure
$\Delta$.
\end{enumerate}
\end{definition}

\noindent Let $\Aut^\ub(\Delta)=\Aut^\ub(\cS,D,\omega)$ and
$\Aut^\ub(\Xi)=\Aut^\ub(\cS,D,J,\omega)$ denote the groups of symmetries of
$\Delta$ and $\Xi$. We have:
\beqa
&& \Aut^\ub(\Xi)=\Aut^\ub(\Delta)\cap \Aut(\cS,J)=\Aut^\ub(\cS,D)\cap \Aut^\ub(\cS,J,\omega)\\
&& \Aut^\ub(\Delta)=\Aut^\ub(\cS,\omega)\cap \Aut^\ub(\cS,D)~~.
\eeqa

\noindent
Given a symplectic automorphism $f\in \Aut^\ub(\cS,\omega)$, the
endomorphism $\bAd(f)(J)$ is again a taming of $(\cS,\omega)$, where
$\bAd(f)$ denotes the adjoint action of $f$ on ordinary sections of
the bundle $End(\cS)$ (see \cite{sigma}). Hence for any
electromagnetic structure $\Xi=(\cS,D,J,\omega)$ having $\Delta$ as
its underlying duality structure, the quadruplet:
\ben
\label{Xif}
\Xi_f\eqdef (\cS,D,\bAd(f)(J),\omega)
\een
is again an electromagnetic structure having $\Delta$ as its
underlying duality structure. This defines a left action of the group
$\Aut^\ub(\cS,\omega)$ on the set $\ES_{\Delta}(\cM)$ of all
electromagnetic structures whose underlying duality structure equals
$\Delta$. 

Let $M$ be a four-manifold and $\cD=(\Sigma,\Xi)$ be a
scalar-electromagnetic structure with underlying scalar structure
$\Sigma=(\cM,\cG,\Phi)$ and underlying electromagnetic structure
$\Xi=(\cS,\cD,J,\omega)$. Let $\cD_0=(\Sigma,\Delta)$ be the
scalar-duality structure underlying $\cD$, where
$\Delta=(\cS,D,\omega)$. Let $g\in \Met_{3,1}(M)$ be a Lorentzian
metric on $M$. Let:
\be
\Aut(\Sigma)\eqdef \{\psi\in \Iso(\cM,\cG)|\Phi\circ \psi=\Phi\}~~,
\ee
where $\Iso(\cM,\cG)$ denote the isometry group of $(\cM,\cG)$. 

\begin{definition} 
The {\em scalar-electromagnetic duality group} of $\cD_0$
is the following subgroup of $\Aut^\ub(\Delta)$:
\be
\Aut(\cD_0) \eqdef \{f\in \Aut^\ub(\Delta)|f_0\in \Aut(\Sigma)\}~~.
\ee
An element of this group is called a {\em scalar-electromagnetic
duality}. The {\em duality action} is the action of $\Aut(\cD_0)$ on
the set $\Conf^g_{\cD_0}(M)$ given by:
\be
\label{dualityaction0}
f \diamond (\varphi, \cV)\eqdef (f_0\circ \varphi, {\hat f}^\varphi (\cV))~~,~~\forall f\in \Aut(\cD_0)~~,
\ee
where $f_0\in \Diff(\cM)$ is the projection of $f$ to $\cM$ and ${\hat f}:\cS\rightarrow \cS^{f_0}$ is 
the based isomorphism of vector bundles induced by $f$. 
\end{definition}

\begin{thm}
\label{thm:DualityInvariance}
For any $f\in \Aut(\cD_0)$, we have:
\ben
f \diamond \Sol^g_{\cD}(M)=\Sol^g_{\cD_f}(M)~~,
\een
where:
\be
\cD_f\eqdef (\Sigma,\Xi_f)~~.
\ee
\end{thm}

\begin{definition}
The {\em scalar-electromagnetic symmetry group} of $\cD$ is the following
subgroup of $\Aut(\cD_0)$:
\be
\Aut(\cD)\eqdef\{f\in \Aut(\cD_0)|\bAd(f)(J)=J\}=\{f\in \Aut^\ub(\Xi)|f_0\in \Aut(\Sigma)\}\, .
\ee
An element of this group is called a {\em scalar-electromagnetic symmetry}.
\end{definition}

\begin{cor}
For all $f\in \Aut(\cD)$, we have:
\be
f\diamond \Sol^g_{\cD}(M)=\Sol^g_{\cD}(M)~~.
\ee
Thus $\Aut(\cD)$ consists of symmetries of the scalar-electromagnetic
equations \eqref{sc} and \eqref{em}, for any fixed
Lorentzian metric $g\in \Met_{3,1}(M)$.
\end{cor}

\noindent We have short exact sequences:
\beqa
&& 1\rightarrow \Aut(\Delta)\hookrightarrow \Aut(\cD_0)\longrightarrow \Aut^{\Delta}(\Sigma)\rightarrow 1~~\\
&& 1\rightarrow \Aut(\Xi)\hookrightarrow \Aut(\cD)\longrightarrow \Aut^{\Xi}(\Sigma)\rightarrow 1~~,
\eeqa
where $\Aut(\Delta)\eqdef \Hom_{\DS(N)^\times}(\Delta,\Delta)$ and
$\Aut(\Xi)\eqdef \Hom_{\ES(N)^\times}(\Xi,\Xi)$ are the groups of {\em
based} automorphisms of $\Delta$ and $\Xi$ and the groups
$\Aut^\Delta(\Sigma)$ and $\Aut^\Xi(\Sigma)$ consist of those
automorphisms of the scalar structure $\Sigma$ which admit lifts to
scalar-electromagnetic duality transformations and
scalar-electromagnetic symmetries, respectively. Let $\Hol_D^p$ be the
holonomy group of $D$ at a point $p\in\cM$. Then we can identify
$\Aut(\Delta)$ with the commutant of $\Hol_D^p$ inside the group
$\Sp(\cS_p,\omega_p)\simeq \Sp(2n,\R)$.

\section{Twisted Dirac quantization}

\noindent Let $N$ be a manifold. 

\subsection{Integral duality structures and integral electromagnetic structures}

\begin{definition}
Let $\Delta=(\cS,D,\omega)$ be a duality structure of rank $2n$
defined on $N$. A {\em Dirac system} for $\Delta$ is a fiber
sub-bundle $\Lambda\subset \cS$ which satisfies the following
conditions:
\begin{enumerate}[1.]
\itemsep 0,0em
\item For any $x\in X$, the triple $(\cS_x,\omega_x,\Lambda_x)$ is an
  integral symplectic space, i.e. $\Lambda_x$ is a full lattice in
  $\cS_x$ and $\omega_x(\Lambda_x,\Lambda_x)\subset \Z$.
\item $\Lambda$ is invariant under the parallel transport of $D$, i.e. 
the following condition is satisfied for any path $\gamma\in \cP(N)$:
\ben
\label{LambdaFlatness}
T_\gamma^\Delta(\Lambda_{\gamma(0)})=\Lambda_{\gamma(1)}~~.
\een
\end{enumerate}
For every $x\in N$, the lattice $\Lambda_x\subset \cS_x$ is called the
{\em Dirac lattice defined by $\Lambda$ at the point $x$}.
\end{definition}

\begin{definition}
An {\em integral duality structure} defined on $N$ is a pair
$\bDelta\eqdef (\Delta,\Lambda)$, where $\Delta$ is is a duality
structure defined on $N$ and $\Lambda$ is a Dirac system for $\Delta$. 
\end{definition}

\noindent Relation \eqref{LambdaFlatness} implies
that the type $\bt$ (the ordered list of elementary divisors) of the
integral symplectic space $(\cS_x,\omega_x,\Lambda_x)$ does not depend
on the point $x\in N$. This quantity is denoted $\bt(\bDelta)$ and
called the {\em type} of $\bDelta$.

\begin{definition}
Let $\bDelta=(\Delta_1,\Lambda_1)$ and
$\bDelta_2=(\Delta_2,\Lambda_2)$ be two integral duality structures
defined on $N$. An morphism of of integral duality structures from
$\bDelta_1$ to $\bDelta_2$ is a morphism of duality structures
$f:\Delta_1\rightarrow \Delta_2$ such that $f(\Lambda_1)\subset
\Lambda_2$.
\end{definition}

\begin{remark}
The set of isomorphism classes of integral duality structures of type
$\bt$ defined on $N$ is in bijection with the character variety:
\be
C_{\pi_1(N)}(\Sp_\bt(2n,\Z))=\Hom(\pi_1(N),\Sp_\bt(2n,\Z))/\Sp_\bt(2n,\Z)~~,
\ee
where $\Sp_\bt(2n,\Z)$ is the modified Siegel modular group of type $\bt$. 
\end{remark}

\noindent Let $\bDelta\eqdef (\cS,D,\omega,\Lambda)$ be an integral duality
structure or rank $2n$ and type $\bt$, defined on $N$.  For any $x\in
N$, the integral symplectic space $(\cS_x,\omega_x,\Lambda_x)$ defines
a symplectic torus $X_s(\cS_x,\omega_x,\Lambda_x)$. These tori fit
into a fiber bundle $\cX_s(\bDelta)$, endowed with a complete flat
Ehresmann connection $\cH_{\bDelta}$ induced by $D$. The Ehresmann
transport of this connection is through isomorphisms of
symplectic tori so it preserves the group structure and
symplectic form of the fibers; in particular, the holonomy group of
$\cH_{\bDelta}$ is contained in $\Sp_\bt(2n,\Z)$.

\begin{definition}
The pair $(\cX_s({\bDelta}),\cH_{\bDelta})$ is called the {\em flat
bundle of symplectic tori} defined by the integral duality structure
$\bDelta$.
\end{definition}

\begin{definition} 
An {\em integral electromagnetic structure} defined on $N$ is a pair
$\bXi=(\Xi,\Lambda)$, where $\Xi=(\Delta,J)$ is an electromagnetic
structure defined on $N$ and $\Lambda$ is a Dirac system for the
underlying duality structure $\Delta=(\cS,D,\omega)$ of $\Xi$.
The type of the integral duality structure $\bDelta=(\cS,D,\omega,
\Lambda)$ is called the {\em type} $\bt(\bXi)$ of $\bXi$:
\be
\bt(\bXi)\eqdef \bt(\bDelta)~~.
\ee
\end{definition}

\noindent Let $\bXi=(\cS,D, J,\omega,\Lambda)$ be an integral
electromagnetic structure of real rank $2n$ and type $\bt$, with
underlying duality structure $\Delta=(\cS,D,\omega)$. For every $x\in
N$, the fiber $(\cS_x,J_x,\omega_x,\Lambda_x)$ is an integral
tamed symplectic space which defines a polarized Abelian variety
$X_h(\cS_x,J_x,\omega_x,\Lambda_x)$ of type $\bt$, whose underlying
symplectic torus is given by $X_s(\cS_x,\omega_x,\Lambda_x)$. These
polarized Abelian varieties fit into a smooth fiber bundle
$\bcX_h(\bXi)$. As above, the connection $D$ induces a complete
integrable Ehresmann connection $\cH_\bXi\eqdef \cH_{\bDelta}$ on this
bundle, whose transport proceeds through isomorphisms of symplectic
tori, so it preserves the Abelian group structure and symplectic form
of the fibers but not their complex structure.

\begin{definition}
The pair $(\cX_h({\bXi}),\cH_\bXi)$ is called the {\em bundle of
  polarized Abelian varieties} defined by the integral electromagnetic
structure $\bXi$.
\end{definition}

\subsection{The twisted Dirac quantization condition}

Let $(M,g)$ be a Lorentzian four-manifold and $(\cM,\cG)$ be a Riemannian manifold. 
Let $\varphi\in \cC^\infty(M,\cM)$ be a fixed
smooth map from $M$ to $\cM$. Let $\bXi=(\Xi,\Lambda)$ be an integral
electromagnetic structure defined on $\cM$, with underlying
electromagnetic structure $\Xi=(\cS,D,J,\omega)$ and underlying
duality structure $\Delta=(\cS,D,\omega)$. Then the system
$\bXi^\varphi=(\Xi^\varphi,\Lambda^\varphi)$ is an electromagnetic
structure on $M$, where $\Lambda^\varphi$ is the $\varphi$-pullback of
the fiber sub-bundle $\Lambda\subset \cS$; this has underlying duality structure
$\Delta^\varphi=(\cS^\varphi,D^\varphi,\omega^\varphi)$. Let
$\bDelta^\varphi\eqdef (\Delta^\varphi,\Lambda^\varphi)$ denote the
integral duality structure underlying $\Xi^\varphi$. Let $\Symp_0$ denote the category of 
finite-dimensional integral symplectic vector spaces. Let
$H^\bullet(M,\bDelta^\varphi)$ denote the total twisted singular
cohomology group of $M$ with coefficients in the
$\Symp_0^\times$-valued local system $T_{\bDelta^\varphi}$ and
let $H^\bullet(M,\Delta^\varphi)$ denote the total twisted singular
cohomology space of $M$ with coefficients in the $\Symp^\times$-valued
local system $T_{\Delta^\varphi}$. The latter can be identified with
the total cohomology space
$H^\bullet_{\dd_{D^{\varphi}}}(M,\cS^\varphi)$ of the twisted de Rham
complex $(\Omega^\bullet(M,\cS^\varphi),\dd_{D^\varphi})$. Since
$\cS^\varphi=\Lambda^\varphi\otimes_\Z\R$, the coefficient sequence
gives a map $j_\ast:H^\bullet(M,\bDelta^\varphi)\rightarrow
H^\bullet(M,\Delta^\varphi)$, whose image
$H^\bullet_{\Lambda^\varphi}(M,\Delta^\varphi)\eqdef
j_\ast(H^\bullet(M,\bDelta^\varphi))$ is a graded subgroup of the
graded additive group of $H^\bullet(M,\Delta^\varphi)$.

\begin{definition}
An electromagnetic field $\cV\in \Omega^2(M,\cS^\varphi)$ is called
{\em $\Lambda^\varphi$-integral} if its $D^\varphi$-twisted cohomology
class $[\cV]\in H^2_{\dd_{D^\varphi}}(M,\cS^\varphi)\equiv
H^2(M,\Delta^\varphi)$ belongs to
$H_{\Lambda^\varphi}^2(M,\Delta^\varphi)$:
\ben
\label{integrality}
[\cV]\in H_{\Lambda^\varphi}^2(M,\Delta^\varphi)=j_\ast(H^2(M,\bDelta^\varphi))~~.
\een
\end{definition}

The condition that $\cV$ be $\Lambda^\varphi$-integral is called 
the {\em twisted Dirac quantization condition} defined by the Dirac
structure $\Lambda$. This should be viewed as a condition constraining
semiclassical Abelian gauge field configurations; a mathematical model
for such configuration can be given using a certain version of twisted
differential cohomology.

\subsection{Integral scalar-electromagnetic duality and symmetry groups}

\begin{definition}
An {\em integral scalar-duality structure} is a pair $\bcD_0\eqdef
(\cD_0,\Lambda)$, where $\cD_0=(\Sigma,\Delta)$ is a scalar-duality
structure and $\Lambda$ is a Dirac system for $\Delta$.
An {\em integral scalar-electromagnetic structure} is a pair
$\bcD=(\cD,\Lambda)$, where $\cD=(\Sigma,\Xi)$ is a
scalar-electromagnetic structure and $\Lambda$ is a Dirac system for
the underlying duality structure of the electromagnetic structure
$\Xi$.
\end{definition}

Let $\bcD=(\cD,\Lambda)$ be an integral scalar-electromagnetic
structure with underlying scalar-electromagnetic structure
$\cD=(\Sigma,\Xi)$, where $\Sigma=(\cM,\cG,\Phi)$ and $\Xi\eqdef
(\cS,D,J,\omega)$. Let $\Delta=(\cS,D,\omega)$ be the underlying
duality structure and let $\bDelta=(\Delta, \Lambda)$ and
$\bXi=(\Xi,\Lambda)$ be the underlying integral duality structure and
integral electromagnetic structure. Let $\cD_0=(\Sigma, \Delta)$ be
the underlying scalar-duality structure and $\bcD_0=(\cD_0,\Lambda)$
be the underlying integral scalar-duality structure.

\begin{definition}
The {\em integral scalar-electromagnetic duality group} defined by the
integral scalar-duality structure $\bcD_0$ is the following subgroup
of the scalar-electromagnetic duality group $\Aut(\cD_0)$:
\be
\Aut(\bcD_0)\eqdef \{f\in \Aut(\cD_0)|f(\Lambda)=\Lambda\}\subset \Aut(\cD_0)~~.
\ee
Elements of this group are called {\em integral scalar-electromagnetic dualities}.
The {\em integral scalar-electromagnetic symmetry group} defined by the
integral scalar-electromagnetic structure $\bcD$ is the following
subgroup of the scalar-electromagnetic symmetry group $\Aut(\cD)$:
\be
\Aut(\bcD)\eqdef \{f\in \Aut(\cD)|f(\Lambda)=\Lambda \}\subset \Aut(\cD)~~.
\ee
Elements of this group are called {\em integral scalar-electromagnetic
  symmetries}.
\end{definition}

\noindent Notice that $\Aut(\bcD)$ is a subgroup of $\Aut(\bcD_0)$. 

\

\noindent{\bf Acknowledgments}
The authors thank Tomas Ortin for discussions and correspondence. The
work of C. I. L. is supported by grant IBS-R003-S1. The work of
C.S.S. is supported by the ERC Starting Grant 259133 -- Observable
String.



\begin{thebibliography}{100}
\bibitem{Ortin}{T.~Ortin, {\em Gravity and Strings}, Cambridge Monographs on Mathematical Physics, 2nd edition, 2015.}
\bibitem{FreedmanProeyen}{D.~Z.~Freedman, A.~Van Proeyen, {\em
Supergravity}, Cambridge Monographs on Mathematical Physics, 2012.}
\bibitem{sigma}{C.~I.~Lazaroiu, C.~S.~Shahbazi, {\em Generalized Einstein-Scalar-Maxwell theories and locally geometric U-folds},  \arxiv{1609.05872}.}
\bibitem{GeometricUfolds}{C. I. Lazaroiu, C. S. Shahbazi, {\em Geometric U-folds in four dimensions}, \arxiv{1603.03095}}
\bibitem{Lipschitz}{C.~I.~Lazaroiu, C.~S.~Shahbazi, {\em Real pinor bundles and real Lipschitz structures}, \arxiv{1606.07894}.}
\bibitem{lip}{C.~I.~Lazaroiu, C.~S.~Shahbazi, {\em On the spin geometry of supergravity and string theory}, \arxiv{1607.02103}.}
\bibitem{BairdWood}{P.~Baird, J.~C.~Wood, {\em Harmonic morphisms between Riemannian manifolds}, Clarendon Press, Oxford, 2003.}
\end{thebibliography}

\end{document}